\documentclass[10pt]{IEEEtran}
\usepackage[utf8]{inputenc}
\usepackage{amssymb}
\usepackage{amsmath}
\usepackage{graphicx}
\usepackage{nomencl}
\usepackage{etoolbox}
\usepackage{algorithm} 
\usepackage{algpseudocode} 
\usepackage{nomencl}
\usepackage[T1]{fontenc}
\usepackage{cite}
\usepackage{array,multirow,graphicx}
\usepackage{multirow}
\usepackage[dvipsnames]{xcolor}
\usepackage{xcolor}

\makenomenclature
\usepackage{url}

\ifCLASSOPTIONcompsoc
  \usepackage[caption=false,font=normalsize,labelfont=sf,textfont=sf]{subfig}
\else
  \usepackage[caption=false,font=footnotesize]{subfig}
\fi

\title{Accelerating L-shaped Two-stage Stochastic SCUC with Learning Integrated Benders Decomposition}
\author{Fouad Hasan, \textit{Student Member, IEEE}, Amin Kargarian, \textit{Senior Member, IEEE}
\date{June 2023 2021}
\thanks{This work was supported by the National Science Foundation under Grant ECCS-1944752}
\thanks{Fouad hasan and A. Kargarian are with the Department of Electrical and Computer Engineering, Louisiana State University, Baton Rouge, LA 70803 USA (e-mail: fhasan1@lsu.edu, kargarian@lsu.edu).}}

\begin{document}
\maketitle
\begin{abstract}
Benders decomposition is widely used to solve large mixed-integer problems. This paper takes advantage of machine learning and proposes enhanced variants of Benders decomposition for solving two-stage stochastic security-constrained unit commitment (SCUC). The problem is decomposed into a master problem and subproblems corresponding to a load scenario. The goal is to reduce the computational costs and memory usage of Benders decomposition by creating tighter cuts and reducing the size of the master problem. Three approaches are proposed, namely regression Benders, classification Benders, and regression-classification Benders. A regressor reads load profile scenarios and predicts subproblem objective function proxy variables to form tighter cuts for the master problem. A criterion is defined to measure the level of usefulness of cuts with respect to their contribution to lower bound improvement. Useful cuts that contain the necessary information to form the feasible region are identified with and without a classification learner. Useful cuts are iteratively added to the master problem, and non-useful cuts are discarded to reduce the computational burden of each Benders iteration. Simulation studies on multiple test systems show the effectiveness of the proposed learning-aided Benders decomposition for solving two-stage SCUC as compared to conventional multi-cut Benders decomposition.
\end{abstract}


\renewcommand\nomgroup[1]{%
  \item[
  \ifstrequal{#1}{I}{\emph{Indices, Sets, and Parameters:}}{%
  \ifstrequal{#1}{V}{\emph{Variables:}}{%
  \ifstrequal{#1}{U}{\emph{Uncertainty Related Parameters and Functions:}}{}}}%
]}

\nomenclature[I]{\(g,t,w,c\)}{Indices for generators, time, stochastic scenarios, and contingency}
\nomenclature[I]{\(k\)}{Benders iteration index}
\nomenclature[I]{\(p_{gtwc}\)}{Generator output}
\nomenclature[I]{\(D_t\)}{Nodal power demand}
\nomenclature[I]{\(f\{\cdot\}\)}{Generation cost function}
\nomenclature[I]{\(SU_g, SD_g\)}{Startup and shutdown costs of unit $g$}

\nomenclature[I]{\(u_{gt}\)}{ON/OFF status of generating unit $g$}
\nomenclature[I]{\(y_{gt}, z_{gt}\)}{Startup and shutdown binary variable indicators}
\nomenclature[I]{\(J_{MP}\)}{Master problem objective function}
\nomenclature[I]{\(J_{SP}\)}{Subproblem objective function}
\nomenclature[I]{\(\alpha_{w}\)}{Proxy of subproblem $ \omega $ objective function}

\nomenclature[I]{\(b\)}{Accumulated cuts}
\nomenclature[I]{\(\gamma_{tw},\nu_{gtw}\)}{Auxiliary continuous variables}
\nomenclature[I]{\(\delta\)}{Cut filtering criterion}
\nomenclature[I]{\(\epsilon\)}{Duality gap tolerance limit}
\nomenclature[I]{\(\eta\{\cdot\}\)}{Random parameter}
\nomenclature[I]{\(\lambda\)}{Dual variables}
\nomenclature[I]{\(\pi_{w}\)}{Probability of stochastic scenario $ \omega $}

\nomenclature[I]{\(\phi\)}{List of filtered cuts}
\nomenclature[I]{\(\psi\{\cdot\}\)}{Numerical value of a cut}

\printnomenclature[0.69in]
\mbox{}



\section{Introduction}
\subsection{Background and Motivation}
\IEEEPARstart{S}{ecurity}-constrained unit commitment (SCUC) is a complex optimization problem solved daily to determine the optimal schedule of generating units \cite{conejo2018unit}. SCUC is a mixed-integer program (MIP) that involves a mix of continuous and integer decision variables. The challenge with mixed integer problems is that the inclusion of integer variables often makes the problem NP-hard. SCUC becomes even more complex under input parameter (e.g., demand) uncertainty \cite{7339456, 9027844, 9899749}  To counterbalance uncertainty in power system operation, SCUC is often formulated as a stochastic MIP problem, known as an L-shaped optimization problem\cite{caroe1998shaped}. An L-shaped problem is formulated as two-stage stochastic optimization. This formulation involves here-and-now and wait-and-see decision variables. It is a two-stage optimization problem with a single integer recourse decision in the second stage. The name L-shaped comes from the shape of the SCUC feasible region when represented graphically. This type of problem is commonly encountered in decision-making under uncertainty, where a decision-maker may face different scenarios or outcomes with associated probabilities. Due to the block structure, the L-shaped two-stage stochastic program is an effective method for SCUC \cite{conejo2006decomposition}. 
L-shaped optimization problems can be solved using specialized techniques such as Benders decomposition (BD). Although Benders decomposition solves L-shaped problems, it suffers from exponential worst-case computational complexity. This is due, in part, to the large number of successively added cuts over iterations, which causes the size of Benders' master problem to grow excessively. Moreover, the master problem takes over 90\% of the time required to implement Benders decomposition \cite{magnanti1981accelerating}. Also, as noted in \cite{minoux1986mathematical}, not every extreme point in the feasible region of subproblems contributes equally to limiting the optimal solution to the master problem. This implies that a significant number of Benders cuts may not be tight enough at the final optimal solution. Consequently, these non-useful Benders cuts can make solving large-scale integer programs challenging.

Various mathematical and heuristic approaches have been presented in the literature to enhance Benders decomposition performance. This study aims to accelerate the convergence of Benders decomposition by taking advantage of machine learning techniques. 
\vspace{-6pt}
\subsection{Literature Review}
It has been more than 50 years since the development of the Benders decomposition algorithm by J. Bender (1962) \cite{bnnobrs1962partitioning}. The algorithm is designed to address complicating variables that, when temporarily fixed, simplify the problem significantly. It distributes the computational load between a master problem and a subproblem (SP). Benders decomposition has proven successful in various fields, including planning and scheduling, healthcare, transportation, telecommunications, energy and resource management, and chemical process design \cite{conejo2006decomposition}. The primary application of Benders decomposition is initially focused on solving MIP problems. Once integer variables are fixed, the problem is converted into a continuous linear program, which can develop cuts using standard duality theory. Many enhancements have been made to extend Benders' applicability to a wider range of problems. As a result, Benders decomposition has been widely used to solve linear, nonlinear, integer, stochastic, multi-stage, and bilevel optimization problems \cite{rahmaniani2017benders}.
The traditional implementation of Benders decomposition can be computationally expensive, time-consuming, and memory-intensive, with problems such as poor feasibility and optimality cuts, ineffective early iterations, and the zigzagging behavior of primal solutions \cite{magnanti1981accelerating}. Researchers have explored various strategies to speed up the convergence and reduce the number of iterations and time spent on each iteration. The master problem is usually solved using branch-and-bound, with the simplex approach used to solve the subproblem. However, a significant number of cuts generated do not contribute to convergence, leading to memory occupation \cite{minoux1986mathematical}. To address this problem, improvement criteria are proposed to ensure that new and useful cuts are included in the master problem \cite{holmberg1990convergence}. Additionally, researchers have observed several orders of improvement when using constraint programming to solve the master problem \cite{correa2007scheduling}. Column generation has been introduced to handle specific structures more effectively and achieve tighter constraints at the root node of the branch-and-bound tree \cite{cordeau2001benders, restrepo2015grammar}. For large subproblems, decomposition, parallelization, and column generation have been used to reduce overall solution time. By adding valid inequalities to the master, one can significantly reduce the number of generated cuts and solution time \cite{naoum2010nested}. Moreover, clustering subproblems can decrease the number of iterations \cite{brandes2011implementierung}.
Recently, researchers have been exploring the use of machine learning (ML) to alleviate the computational complexity of complex problems \cite{paulus2021comboptnet, mohammadi2023empowering}. ML-based algorithms can be categorized into three types: end-to-end learning (with label), unsupervised/reinforcement learning (no ground truth), and algorithm-specific hybrid learning, which leverages the specific structure of the target problem to accelerate the optimization process. For instance, a support vector machine is used in \cite{jia2021benders} to construct a cut classifier that identifies valuable cuts in each Benders iteration, thus reducing the size of the master problem and shortening the solution time. Another strategy is to use the Lagrange multipliers of the Benders subproblem to aggregate the optimality cuts \cite{vandenbussche2019data}. In \cite{lee2020accelerating}, a cut regressor and a cut classifier are used, and five features are designed to characterize the generated cuts \cite{ruszczynski2003decomposition}.
Solving two-stage stochastic SCUC with Benders decomposition is still computationally challenging and needs further investigation. Also, despite machine learning advantages, the literature lack ML-based Benders approaches to accelerate two-stage stochastic SCUC.
\vspace{-9pt}
\subsection{Contribution}
This paper presents combined machine learning and model-based variants of Benders decomposition for solving two-stage stochastic SCUC. The proposed approaches form tighter cuts and a truncated, less computationally expensive Benders master problem. A regressor is trained to predict proxy variables corresponding to subproblems and create a tighter master problem feasible region. After carrying out each Benders iteration, useful cuts are identified using a criterion and a classification learner. A truncated master problem is created by adding useful cuts to the model and discarding non-useful cuts. Combining the regressor and classifier reduces the computational cost and memory usage of the Benders decomposition for solving two-stage stochastic SCUC. The Benders iterative loop inherently ensures the feasibility and optimality of results. The three enhanced Benders approaches are tested on various test systems, and the results are analyzed.

\subsection{Paper Organization}
The rest of the paper is structured as follows. The problem formulation is given in Section II. The proposed enhanced Benders approaches are presented in Section III. The numerical simulations are discussed in Section IV, and concluding remarks are provided in Section V.

\section{TWO-STAGE STOCHASTIC UNIT COMMITMENT}
The considered problem is a two-stage stochastic unit commitment with respect to demand uncertainty. Benders decomposition is used to solve the problem.

\subsection{Problem Formulation}
The problem is formulated in (1), which includes two sets of variables pertaining to here-and-now decisions and wait-and-see decisions \cite{conejo2018power}. On/off status of units are here-and-now variables, and generation dispatches are wait-and-see variables. The first term of (1a) is the first-stage startup and shutdown costs, and the second is the second-stage generation dispatch costs. The first stage constraints are (1b) – (1i), which are scenario-independent unit commitment constraints. The second-stage scenario-dependent operational constraints are formulated in (1j)–(1n). N-1 security constraints are included in the model with index c.

\begin{equation} \tag{1a}
    \label{1a}
    \min  \sum_t \sum_g (SU_g y_{gt} + SD_g z_{gt}) + \sum_\omega \pi_\omega \sum_t\sum_g f(p_{gt\omega}) 
\end{equation}
s.t. 
\begin{equation}\tag{1b}
    \label{1b}
    y_{gt}-z_{gt}=u_{gt}-u_{g(t-1)}   \ \ \forall g, \forall t
\end{equation}
\begin{equation}\tag{1c}
    \label{1c}
    y_{gt}+z_{gt} \leq 1 \ \ \forall g, \forall t
\end{equation}
\begin{equation}\tag{1d}
    \label{1d}
    \sum_{t=1}^{UT_g} (1-u_{gt})=0 \ \ \forall g
\end{equation}
\begin{equation}\tag{1e}
    \label{1e}
    \sum_{\tau=t}^{t+T_g^{on}-1} u_{g\tau} \geq T_g^{on}y_{gt} \ \ \forall g, t=UT_g+1, UT_g+2,... T-T_g^{on}+1 
\end{equation}
\begin{equation}\tag{1f}
    \label{1f}
    \sum_{\tau=t}^{T} (u_{g\tau}-y_{gt}) \geq 0 \ \ \forall g, t=T-T_g^{on}+2,...T
\end{equation}

\begin{equation}\tag{1g}
    \label{1g}
    \sum_{t=1}^{DT_g} u_{gt}=0 \ \ \forall g
\end{equation}
\begin{equation}\tag{1h}
    \label{1h}
    \sum_{\tau=t}^{t+T_g^{off}-1} (1-u_{g\tau}) \geq T_g^{off}z_{gt} \ \ \forall g, t=DT_g+1, DT_g+2,...T-T_g^{off}+1 
\end{equation}
\begin{equation}\tag{1i}
    \label{1i}
    \sum_{\tau=t}^{T} (1-u_{g\tau}-z_{gt}) \geq 0 \ \ \forall g, t=T-T_g^{off}+2,...T
\end{equation}

\{
\begin{equation}\tag{1j}
    \label{1j}
    p_g^{min}u_{gt}\leq p_{gt\omega c} \leq p_g^{max}u_{gt} \ \ \forall g, \forall t
\end{equation}
\begin{equation}\tag{1k}
    \label{1k}
    p_{gt\omega c} - p_{g(t-1)\omega c} \leq RU_g(1-y_{gt}) + P_g^{min}y_{gt} \ \ \forall g, \forall t 
\end{equation}
\begin{equation}\tag{1l}
    \label{1l}
    p_{g(t-1)\omega c} - p_{gt\omega c}  \leq RD_g(1-z_{gt}) + P_g^{min}z_{gt} \ \ \forall g, \forall t 
\end{equation}\begin{equation}\tag{1m}
    \label{1m}
    \sum_g p_{gt\omega c} = D_t \ \ \forall t
\end{equation}
\begin{equation}\tag{1n}
    \label{1n}
    - PL^{max} \leq SF_c P_{t\omega c}^{inj} \leq PL^{max} \ \ \forall t
\end{equation}
\}

\begin{equation} 
    \forall \omega \in \Omega, \forall c  \nonumber
\end{equation}

\subsection{Benders Decomposition }
The computational complexity of mixed-integer program (1) increases with the increasing size of the network, number of demand scenarios, and number of contingencies. Benders decomposition is suitable for solving MIP problems with a block structure over contingencies and uncertainty scenarios \cite{conejo2006decomposition}. We use a multi-cut variant of Benders decomposition that converges faster than the classical single-cut approach \cite{vandenbussche2019data, birge1988multicut}. The multi-cut technique decomposes a problem into multiple subproblems and generates multiple cuts in each iteration. This multi-cut generation improves convergence performance \cite{su2015computational}. 
Problem (1) is decomposed into a master problem and n subproblems. The master problem (MP) at iteration k is formulated in (2a)-(2c).

\begin{equation} \tag{2a}
    \label{2a}
    \min_x  \sum_t \sum_g (SU_g y_{gt} + SD_g z_{gt}) + \sum_\omega \alpha_{\omega}
\end{equation}
s.t.
\begin{equation} \tag{2b}
    \label{2b}
    (1b)-(1i)
\end{equation}
\begin{equation} \tag{2c}
    \label{2c}
    \cup_{i=1}^{i=k-2} h_{\beta}^n(u_{gt},\alpha_{\omega})
\end{equation}
\begin{equation} \tag{2d}
    \label{2d}
    h_{\beta}^{k-1}(u_{gt},\alpha_{\omega}): \alpha_\omega  \geq J_{SP,\omega}^{k-1} + \sum_t \sum_g \lambda_{gt\omega}^{k-1} (u_{gt}-u_{gt}^{k-1}) \ \ \forall \omega
\end{equation}
\begin{equation} \tag{2e}
    \label{2e}
    \alpha_\omega \geq \alpha_\omega^{min} \ \ \forall \omega 
\end{equation}
\begin{equation} 
     x \in \{u_{gt},y_{gt},z_{gt},\alpha_\omega\}   \nonumber
\end{equation}

The objective function (2a) consists of startup and shutdown costs and a term as a proxy for the subproblems' objective function. Equation (2c) is the accumulation of all cuts generated up to iteration k-2, (2d) denotes Benders cuts generated at iteration k-1, and (2e) defines the bound on proxy variables where $\alpha_{\omega}^{min}$ is set as a large negative constant. 

Benders subproblem $\omega$ ($SP_\omega$) is formulated in (3a)- (3e).  To gain extra computational efficiency, we use an always-feasible subproblem model \cite{nasri2015network}. This model includes non-negative slack variables $\nu$ and $\gamma$ to prevent ($SP_\omega$) infeasibility. Auxiliary variables $\nu$ and $\gamma$ relax generation capacity constraints.

\begin{equation} \tag{3a}
    \label{3a}
    \min_x  \sum_t \sum_g \sum_\omega (\pi_\omega f(p_{gt\omega})+\nu_(gt\omega)) + \sum_t \sum_\omega \gamma_{t\omega}
\end{equation}
s.t.

\begin{equation} \tag{3b}
    \label{3b}
    u_{gt}=u_{gt}^*
\end{equation}
\{
\begin{equation} \tag{3c}
    \label{3c}
    (1j)-(1n)
\end{equation}
\begin{equation} \tag{3d}
    \label{3d}
    p_g^{min}u_{gt}\leq p_{gt\omega c} \leq p_g^{max}u_{gt} + \nu_{(gt\omega)} + \gamma_{t\omega} \ \ \forall g, \forall t
\end{equation}
\begin{equation} \tag{3e}
    \label{3e}
    \nu_{(gt\omega)} \geq 0, \gamma_{t\omega} \geq 0
\end{equation}
\}

\begin{equation} 
    x \in \{p_{gt\omega c}, \nu_{(gt\omega)}, \gamma_{t\omega}\} \nonumber
\end{equation}

The objective function (3a) consists of generation costs and constraint violation penalty modeled by $\nu_{gtw}$ and $\gamma_{tw}$. Constraint (3b) sets the first stage decisions, i.e., generator unit status, as fixed values received from the master problem. Inequalities (3e) set bound on auxiliary variables. The MP and SPs are solved iteratively until convergence tolerance $ \varepsilon $  is smaller than a predetermined threshold. 

\begin{equation} \tag{4}
    \label{4}
    \epsilon = \frac{UB-LB}{abs(LB)}
\end{equation}
where 
\begin{equation} 
     LB=J_{MP}, \ \ 
     UB=J_{SP}+J_{MP}-\sum_\omega \alpha_\omega \nonumber
\end{equation}

The two-stage stochastic SCUC with multi-cut Benders is summarized in Algorithm I.

\begin{figure}
\centering
    \includegraphics[width=.5\textwidth]{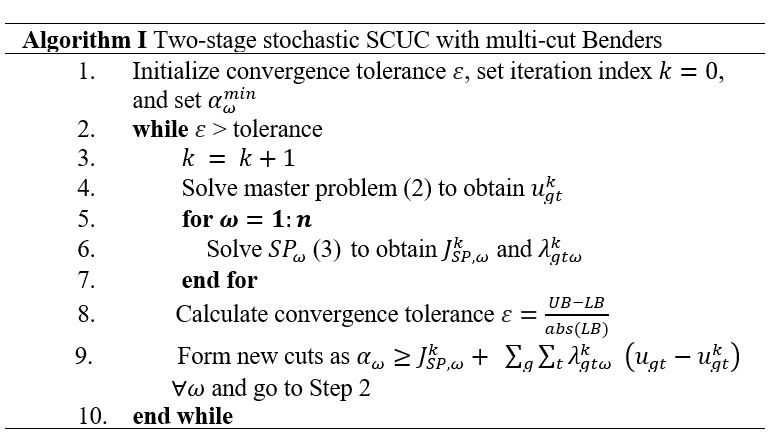}
\label{Algorithm1}
\vspace{-40pt}
\end{figure}

\section{PROPOSED LEARNING ASSISTED BENDERS DECOMPOSITION}
While multi-cut Benders outperforms single-cut Benders in terms of convergence performance, its memory usage and computational overhead still need enhancement. We present three learning-aided approaches to reduce the computational costs of multi-cut Benders decomposition for solving two-stage stochastic SCUC. We mainly focus on the master problem, which has higher computational costs than subproblems.

\subsection{Regression-assisted Benders (R-Benders)}
The first proposed approach uses a combination of a regressor and cut filtering. The regressor forms a tighter bound for the master problem at the root node starting from k=1. We predict the optimal subproblem proxy variables $\alpha_{\omega}^*$, i.e., $\alpha_{\omega}$ upon the convergence of Algorithm I. Knowing $\alpha_{\omega}^*$, or even a good approximation, makes inequalities (2c), (2d), and (2e) tighter and the lower bound LB closer to the optimal subproblem cost. The cut filtering approach identifies and drops non-useful cuts during each subsequent Benders iteration to reduce the size of the master problem. Non-useful cuts are cuts with no impact on the master problem feasible region.

\subsubsection{SP Objective Proxy Prediction}
In the multi-cut Benders formulation (2), $\alpha_{\omega}$ is a variable and $\alpha_{\omega}^{min}$ is a bound where $\alpha_{\omega}$ $\geq$ $\alpha_{\omega}^{min}$ constitutes a master problem constraint. The value of $\alpha_{\omega}$ increases after each iteration and reaches its optimal value $\alpha_{\omega}^*$ upon convergence. The value of $\alpha_{\omega}^{min}$ can be selected by analyzing the physical and economic aspects of the SCUC problem. Setting a suitable $\alpha_{\omega}^{min}$ reduces the search space and enhances Benders convergence performance. An ideal case is to set $\alpha_{\omega}^{min}$ corresponding to each subproblem $ \omega $ as the optimal value of $\alpha_{\omega}$ instead of using a large negative value. Fig. 1 shows the advantage of using optimal $\alpha_{\omega}^*$. The master problem objective reaches its optimal value in the very initial iterations. We use machine learning to predict $\alpha_{\omega}^*$ for each subproblem $\omega$ by reading power demand before implementing multi-cut Benders Algorithm I.

\begin{figure}
\centering
    \includegraphics[width=.3\textwidth]{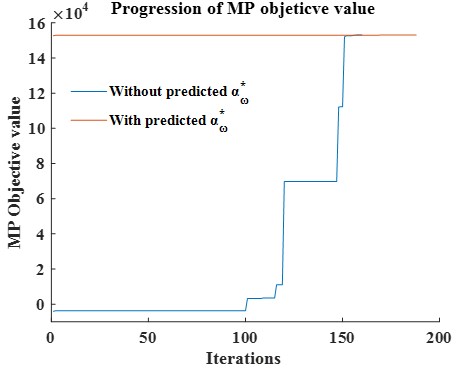}
\caption{Benders convergence using $\alpha_{\omega}^*$  $\forall$ $\omega$ in inequality (2b).}
\label{fig1}
\end{figure}

\subsubsection{Training Dataset Preparation}
Electricity demand scenarios are input to two-stage SCUC. We have observed that $\alpha_{\omega}^*$ are correlated with system demand. We use demand information as the input to a supervised learning model whose output is  $\alpha_{\omega}^*$. To account for potential system operation scenarios in the training phase, we generate a set of daily load profile samples$D^s$ \(\in \mathbb{R}^{s\times 1}\) according to (5). To follow the concept of demand uncertainty modeling in the two-stage SCUC problem, each $D^s$ contains n demand profile scenarios generated by (6). 

\begin{equation} \tag{5}
    \label{5}
    D^s = D_{base}[\eta_s^L+\eta_s(\eta_s^U-\eta_s^L)] \ \ \forall s
\end{equation}
\begin{equation} \tag{6}
    \label{6}
    D^s_\omega = D^s[\eta_\omega^L+\eta_\omega(\eta_\omega^U-\eta_\omega^L)] \ \ \forall \omega, \forall s
\end{equation}

where $D_{\omega}^s \in \mathbb{R}^{s \times 1}$, and \(\eta_{s}\{\cdot\}\)  and \(\eta_{\omega}\{\cdot\}\) follows a uniform distribution between 0 and 1. Indices for samples and subsamples are denoted by $s$ and $ \omega $. Two stages of randomness are carried out to consider a range of realistic operational conditions. The load point is shifted using the random parameter \(\eta_{s}\{\cdot\}\) to model its daily and seasonal volatility within the range [\(\eta_{s}^U\), \(\eta_{s}^L\)]. \(\eta_{s}^U\)  and \(\eta_{s}^L\) can be determined by historical data and load growth projection or to the point when any further increase or decrease renders unit commitment infeasible. The uncertainty in subsamples is modeled using the random parameter \(\eta_{\omega}\{\cdot\}\) within a specified range [\(\eta_{\omega}^U\), \(\eta_{\omega}^L\)]. Subsamples model the hourly load randomness by multiplying $D^s$ by random parameters and creating $D_{\omega}^s$. Fig. 2 illustrates the load scenario generation.

\begin{figure}
\centering
    \includegraphics[width=.5\textwidth]{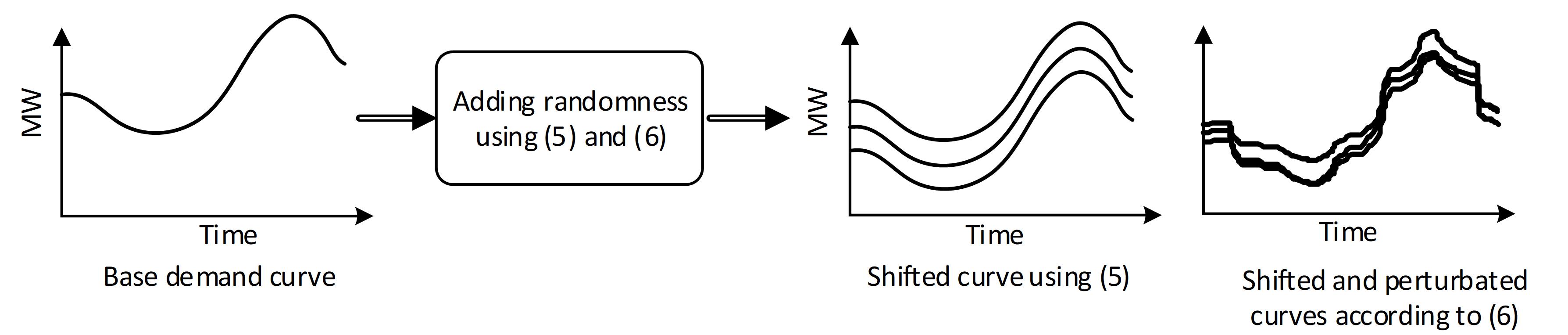}
\caption{Demand scenario generation.}
\label{fig2}
\vspace{-15pt}
\end{figure}

For a given system, Algorithm I is executed for each demand profile sample $D^s$. Optimal values of proxy variables upon Benders convergence are labeled as $\alpha^{s*}$ and stored. Expressions (7) and (8) are, respectively, the input and target of the supervised learning model.

\begin{equation} \tag{7}
    \label{7}
     Input \ \ sample: D^s = [D_1^s, D_2^s, D_3^s,...,D_n^s]^T \ \ \forall s
\end{equation}

\begin{equation} \tag{8}
    \label{8}
    Target \ \ sample: \alpha^{s*}=[\alpha_1^{s*},\alpha_2^{s*},\alpha_3^{s*},...,\alpha_n^{s*}]^T \ \ \forall s  
\end{equation}

\subsubsection{Supervised Learning Strategy}
We use neural networks (NN), an efficient tool to capture the complexity and nonlinearity of a function by utilizing various activation functions. The regressor uses a fully connected NN with mini-batch gradient descent and Rectified Linear Units (ReLU) activation functions for hidden layers. The loss function is the mean squared error (MSE).

\begin{equation} \tag{9}
    \label{9}
     MSE = \frac{\sum_n (x_n-\hat{x_n})^2}{n} 
\end{equation}

The Adam optimizer determines optimal weights and trains the learning model. Various batch sizes, epochs, and layer counts are tested to determine the best architecture. A single hidden layer is selected to have a minimalistic model. Table I shows the architecture and hyperparameters used in this study. Although more sophisticated structures may improve the results, we obtained promising outcomes with this minimalistic architecture. Fig. 3 illustrates a conceptual schematic of the learner model.

\begin{table}
\setlength{\arrayrulewidth}{0.3mm}
\begin{center}
\caption{Hyperparameters of The Learner}
\label{table1}
\begin{tabular}{ m{2cm}  m{2.5cm}  m{3.5cm}} \hline
NN regressor& Hidden layer=1, Batch size=300$~$500  & Activation = ReLU, Loss function = MSE Optimizer = Adam \\ \hline
\end{tabular}
\end{center}
\end{table}

\begin{figure}
\centering
    \includegraphics[width=.3\textwidth]{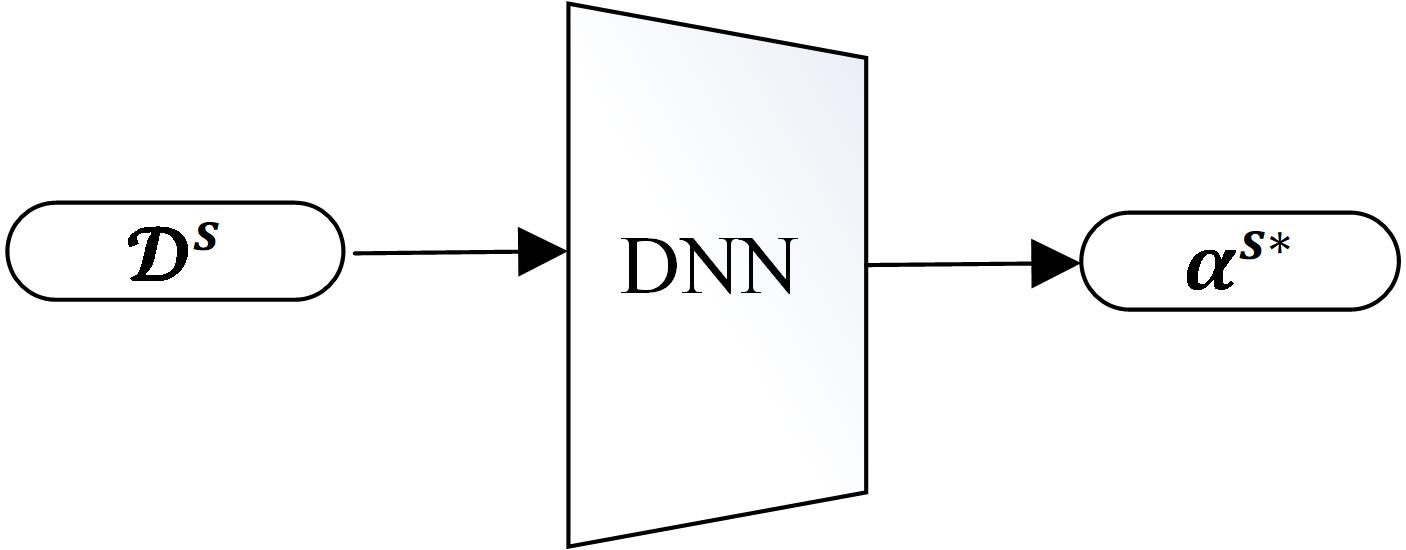}
\caption{NN regressor reads demand profiles and predicts proxy variables $ \alpha $.}
\label{fig3}
\end{figure}

\subsubsection{Data Scaling}
Data normalization and standardization are crucial before training a learner. The data samples are normalized by (10). This process enhances prediction accuracy and numerical stability.

\begin{equation} \tag{10}
    \label{10}
     x_{normalized} = \frac{ x_n-x_{min}}{x_{max}-x_{min}} 
\end{equation}

\subsubsection{Useful Cut Identification}
The master problem still retains all cuts generated after carrying out each Benders iteration. As a result, the master problem size grows at every iteration, calling for a considerable computational memory requirement. In most cases, a subset of Benders cuts generated at each iteration contains the necessary information to build the feasible search space of the master problem. Currently, a lack of practical and systematic approach exists for classifying useful and non-useful cuts for large-scale problems \cite{rahmaniani2017benders}. Several features are suggested in \cite{jia2021benders, lee2020accelerating}, such as cut violation, cut depth, cut order, and cut producing scenario, to identify approximately useful cuts.

Our idea of useful cut identification is based on Fig. 4, which shows an example of the master problem objective value (or lower bound) improvement with cumulatively added cuts. It can be observed that all cuts are not equally contributing to the lower bound and optimality gap improvement. Several cuts provide a positive increase in the lower bound. Such cuts can be classified as useful cuts. Other cuts do not contribute to the lower bound improvement and can be classified as non-useful. Our observation and experiment show that the numerical difference between cut values $\psi{(u_{gt}^k)}$ and proxy values $\alpha^k_\omega$ can be used for useful cut identification. If an inequality constraint is satisfied as equality upon solving optimization, it is typically referred to as a binding or strictly useful constraint. Thus, if $\alpha^k_\omega$ - $\psi{(u_{gt}^k)}$ = 0, the cut is useful. 

\begin{figure}
\centering
    \includegraphics[width=.3\textwidth]{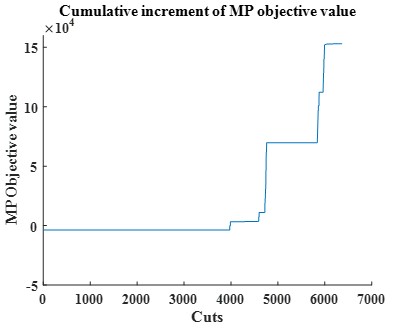}
\caption{Benders lower bound improvement with cuts added cumulatively; a case of IEEE 118-bus system.}
\label{fig4}
\end{figure}

However, many other cuts that are not satisfied as equality may also contain the necessary information to form the master problem feasible region. Such cuts should also be considered useful. As an inequality constraint approaches equality, its importance in forming a feasible region is expected to increase. Consider a demand profile sample $D^s$. To capture all cuts that contain necessary information at each iteration k, we label cuts as "useful" if the difference between $\alpha^k_\omega$ and $\psi{(u_{gt}^k)}$ is less than a threshold $\delta$, which is selected through experimental observation. We have tested various cases and observed that $\delta$ within the range $[0.1,100]$ works well for the studied cases. Once all cuts are generated at iteration $k-1$, we add them to MP and solve the updated MP at iteration k to obtain $u_{gt}^k$ and $\alpha^k_\omega$. We check (11) to identify useful cuts generated at iteration k-1. Cuts that satisfy (11) contain necessary information and are useful. We repeat this process at each iteration $k$ to detect useful cuts generated at iteration $k-1$. Non-useful cuts are dropped from MP before moving to iteration $k+1$. 

\begin{equation} \tag{11}
    \label{11}
    h_{\beta}^{k-1}(u_{gt},\alpha_{\omega})^\phi: |\alpha_\omega^k - \psi(u_{gt}^k)| \leq \delta  
\end{equation}

where $\psi{(u_{gt}^k)}$ is the numerical value of the right-hand side of (2d).

\begin{equation} 
     \psi(u_{gt}^k) = J_{SP,\omega}^{k-1} + \sum_t \sum_g \lambda_{gt\omega}^{k-1} (u_{gt}-u_{gt}^{k-1}) \nonumber
\end{equation}

\begin{figure}
\centering
    \includegraphics[width=.4\textwidth]{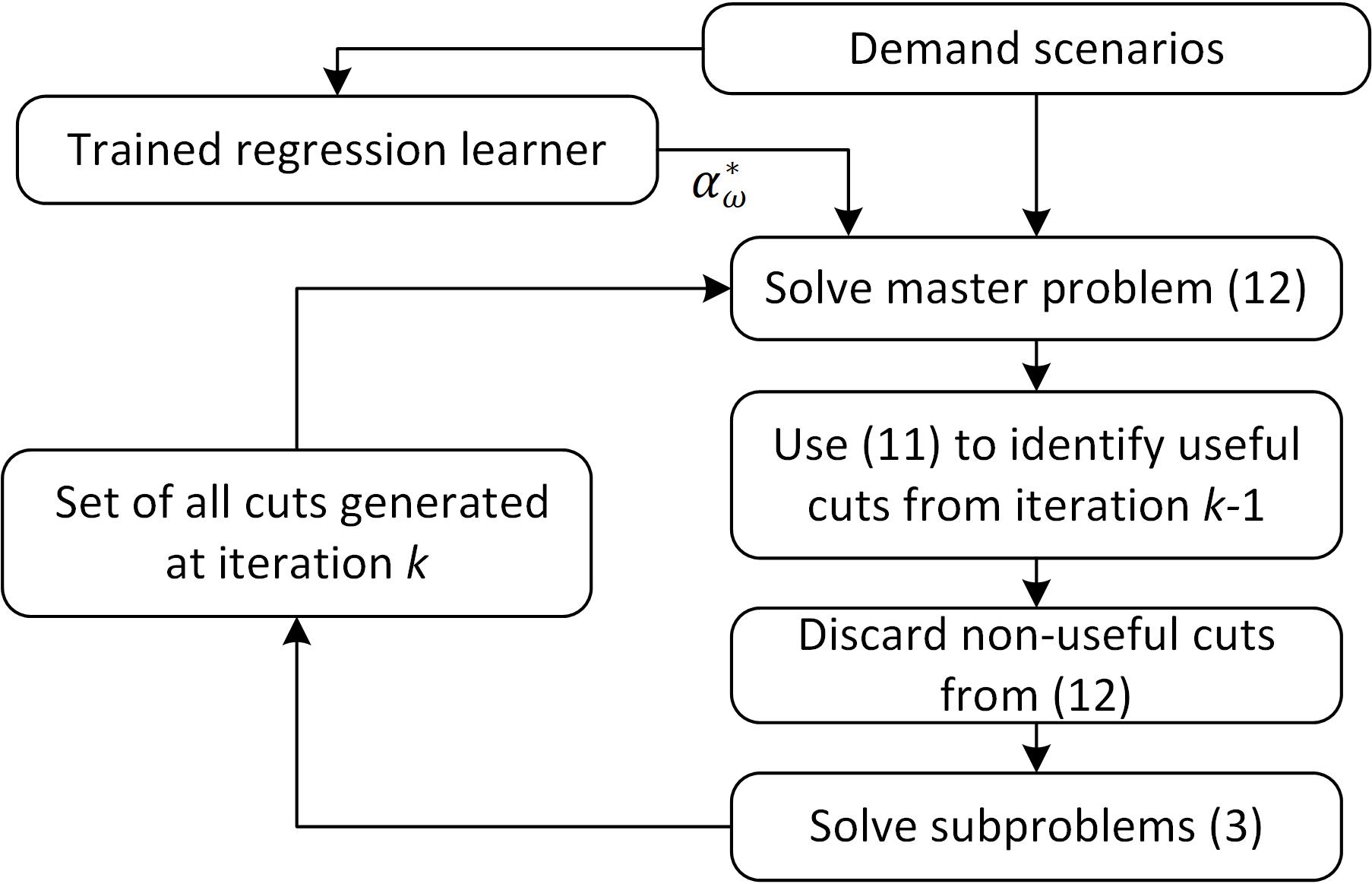}
\caption{R-Benders flowchart.}
\label{fig5}
\end{figure}

\subsubsection{Avoiding Information Loss}
Selecting a small $\delta$ results in labeling fewer cuts as useful. A large $\delta$ yields keeping more cuts in the master problem, and thus less computational cost reduction is gained. Depending on the value of hyperparameter $\delta$, we have observed cases with no cut labeled as useful at a few iterations. Retaining at least one cut from each iteration is crucial; otherwise, the proposed R-Benders may not converge due to loss of information from that particular iteration. To prevent this, we keep several cuts generated from high-load scenarios at each iteration, most of which may already be on the useful cut list. The number of retained cuts is determined through experimentation. 

\subsubsection{R-Benders MP Formulation}
The R-Benders master problem at iteration k is as follows: 

\begin{equation} \tag{12a}
    \label{12a}
    \min_x  \sum_t \sum_g (SU_g y_{gt} + SD_g z_{gt}) + \sum_\omega \alpha_{\omega}
\end{equation}
s.t.
\begin{equation} \tag{12b}
    \label{12b}
    (1b)-(1i)
\end{equation}
\begin{equation} \tag{12c}
    \label{12c}
    \cup_{i=1}^{i=k-2} h_{\beta}^n(u_{gt},\alpha_{\omega})^\psi
\end{equation}
\begin{equation} \tag{12d}
    \label{12d}
    h_{\beta}^{k-1}(u_{gt},\alpha_{\omega}): \alpha_\omega  \geq J_{SP,\omega}^{k-1} + \sum_t \sum_g \lambda_{gt\omega}^{k-1} (u_{gt}-u_{gt}^{k-1}) \ \ \forall \omega
\end{equation}
\begin{equation} \tag{12e}
    \label{12e}
    \alpha_\omega \geq \alpha_\eta \alpha_\omega^{min} \ \ \forall \omega 
\end{equation}
\begin{equation} 
     x \in \{u_{gt},y_{gt},z_{gt},\alpha_\omega\}   \nonumber
\end{equation}

where (12c) denotes the list of all useful cuts detected until iteration $k-2$. (12d) includes all cuts generated at iteration $k-1$. The value of $\alpha_{\omega}$ successively progresses to reach the optimal value. To maintain the quality of the lower bound solution, we should ensure $\alpha_{\omega}$$\geq$$\alpha_{\omega}^*$. Given the possibility of machine learning error, the regressor predictions $\alpha_{\omega}^*$ are reduced by a factor $\alpha_{\eta}$  and $\alpha_{\omega}$$\geq$$\alpha_{\omega}^*$ is replaced with $\alpha_{\omega}$$\geq$$\alpha_{\eta}$$\alpha_{\omega}^*$. The reduction factor $\alpha_{\eta}$ is determined by analyzing the regressor prediction error. 
The R-Benders approach is summarized in Algorithm II and Fig. 5. Although a classification learning model can be trained to identify useful cuts (as in the C-Bender approach presented in Section III.B), the R-Benders approach does not use any classifier with the cost of filtering cuts with one iteration delay. The only learning model in R-Benders is a regressor to predict $\alpha_{\omega}^*$.

\begin{figure}
\centering
    \includegraphics[width=.5\textwidth]{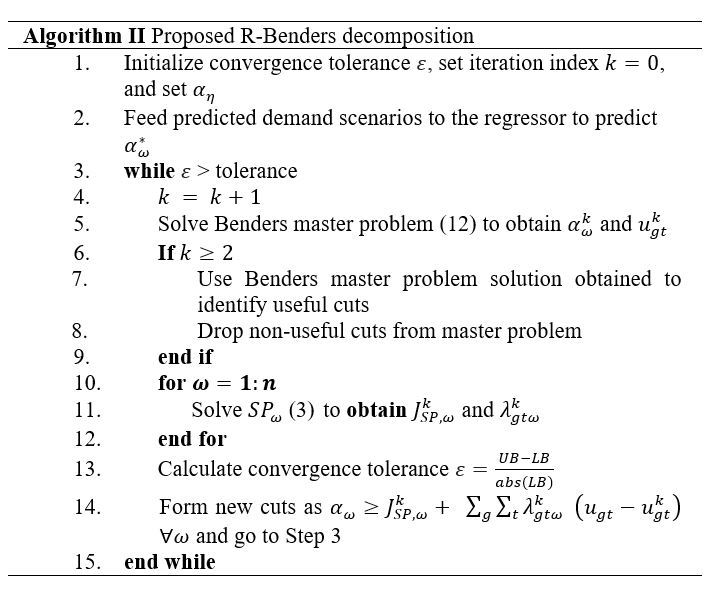}
\label{Algorithm2}
\end{figure}

\subsection{Classification-based Cut Identification (C-Benders)}
We use a supervised classification learner to identify useful cuts before solving the master problem at each iteration k, instead of using criterion (11). Fig. 6 conceptualizes the proposed C-Benders where dotted lines denote non-useful cuts and solid lines denote useful cuts. The input to the classifier is cuts generated at iteration k, and its output is cut labels. Only useful cuts are kept; the rest are omitted from the master problem feasible region. A neural network is trained with a suitable loss function, such as F-score. The model architecture and hyperparameters are given in Table I.

\begin{figure}
\centering
    \includegraphics[width=.4\textwidth]{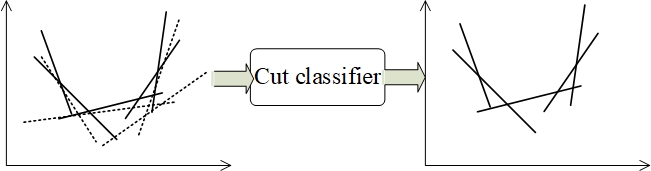}
\caption{Useful cut identification with a classifier.}
\label{fig6}
\end{figure}

To prepare a training dataset, all cuts must be labeled as useful or non-useful. We implement Benders for each demand profile sample and collect all cuts upon convergence. We restart the Benders only to solve the master problem while ignoring subproblems as we have already stored cuts. At each iteration k, we add to the master problem one of the stored cuts and measure its contribution to lower bound improvement. We repeat this procedure by cumulatively adding cuts one after another to label all cuts of iteration k depending on their contribution to the lower bound improvement. We repeat this process for all Benders iterations. One alternative is to form the master problem with all cuts first and then remove cuts one by one to observe their impact on the lower bound.
Generating a labeled dataset for C-Benders is computationally expensive, particularly for multi-cut Benders decomposition. The number of Benders iterations is relatively high for the SCUC problem with intertemporal constraints, leading to a multiplicatively larger number of cuts that should be individually checked. The number of cuts generated at each iteration depends on the number of demand scenarios. This is another factor that increases the computational cost of cut labeling. For instance, consider solving the 118-bus system with a 4-hour commitment horizon and a load profile sample with 40 stochastic scenarios. This case takes 289 iterations to converge with a less than 1\% gap. The total number of cuts generated for each sample is 40×289 =11,560. To determine the label of each cut, we must solve the master problem 11,560 times. Generating the training dataset for C-Benders is computationally expensive, while R-Benders does not need such dataset generation. However, C-Benders provides more accurate cut labels than R-Benders, which estimates labels using (11).
Fig. 7 shows the C-Benders flowchart. The master problem is formulated in (12a)-(12d), and (2e). Filtered cuts in (12c) are identified using the classifier. The training process is offline and needs to be performed only once. The overall accuracy of the learning-aided cut classification is slightly better than the non-learning cut classification of R-Benders. 

\begin{figure}
\centering
    \includegraphics[width=.3\textwidth]{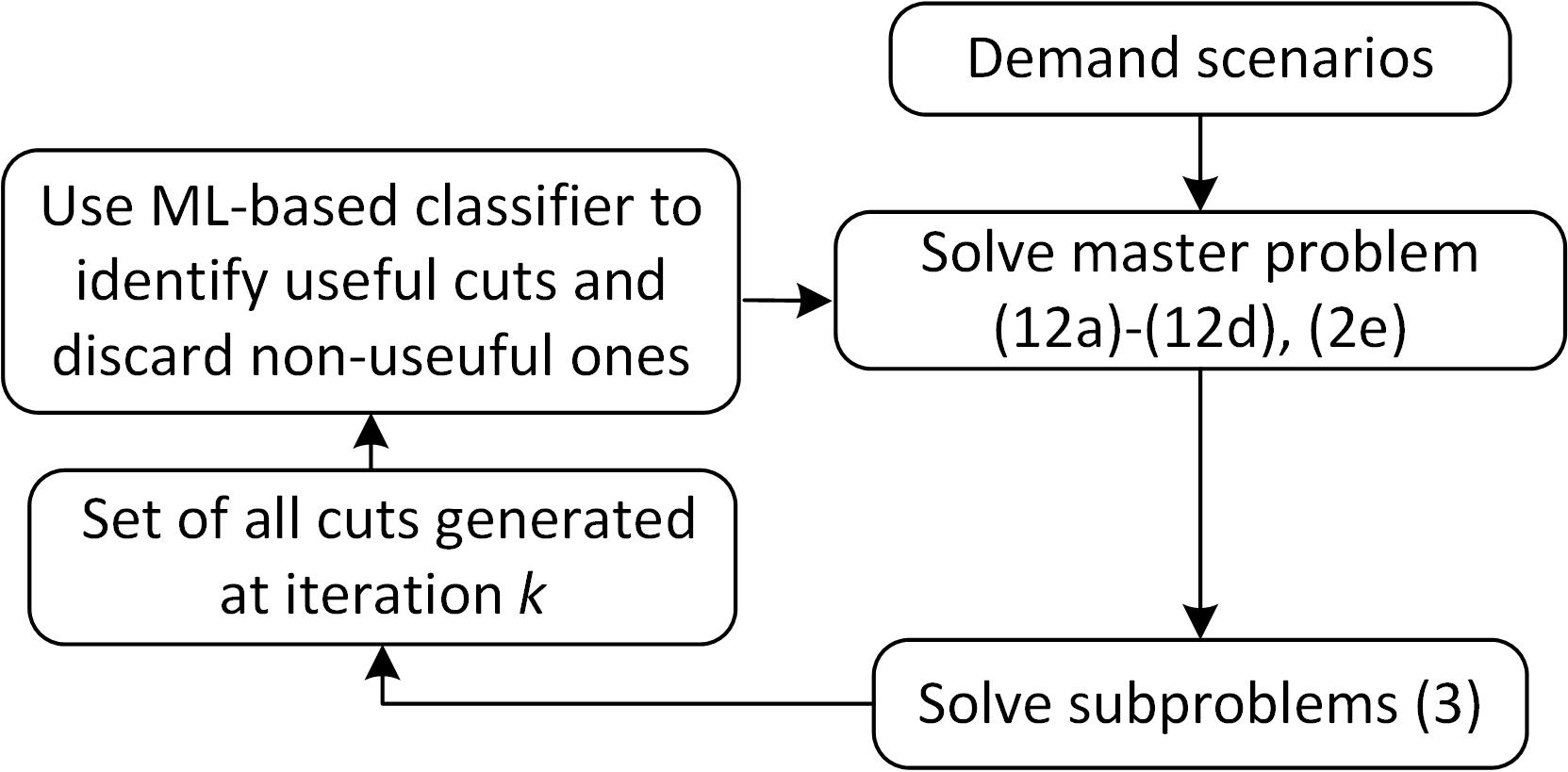}
\caption{C-Benders flowchart.}
\label{fig7}
\end{figure}

\subsection{Combined Classifier-Regressor Benders (CR-Benders)}
This approach (Fig. 8) is a combination of C-Benders and R-Benders. CR-Benders uses a regressor to predict subproblem proxy variables $\alpha_{\omega}^*$ and a classifier for cut filtering before starting each iteration k. This approach is C-Benders equipped with $\alpha_{\omega}^*$ prediction capability.

\begin{figure}
\centering
    \includegraphics[width=.3\textwidth]{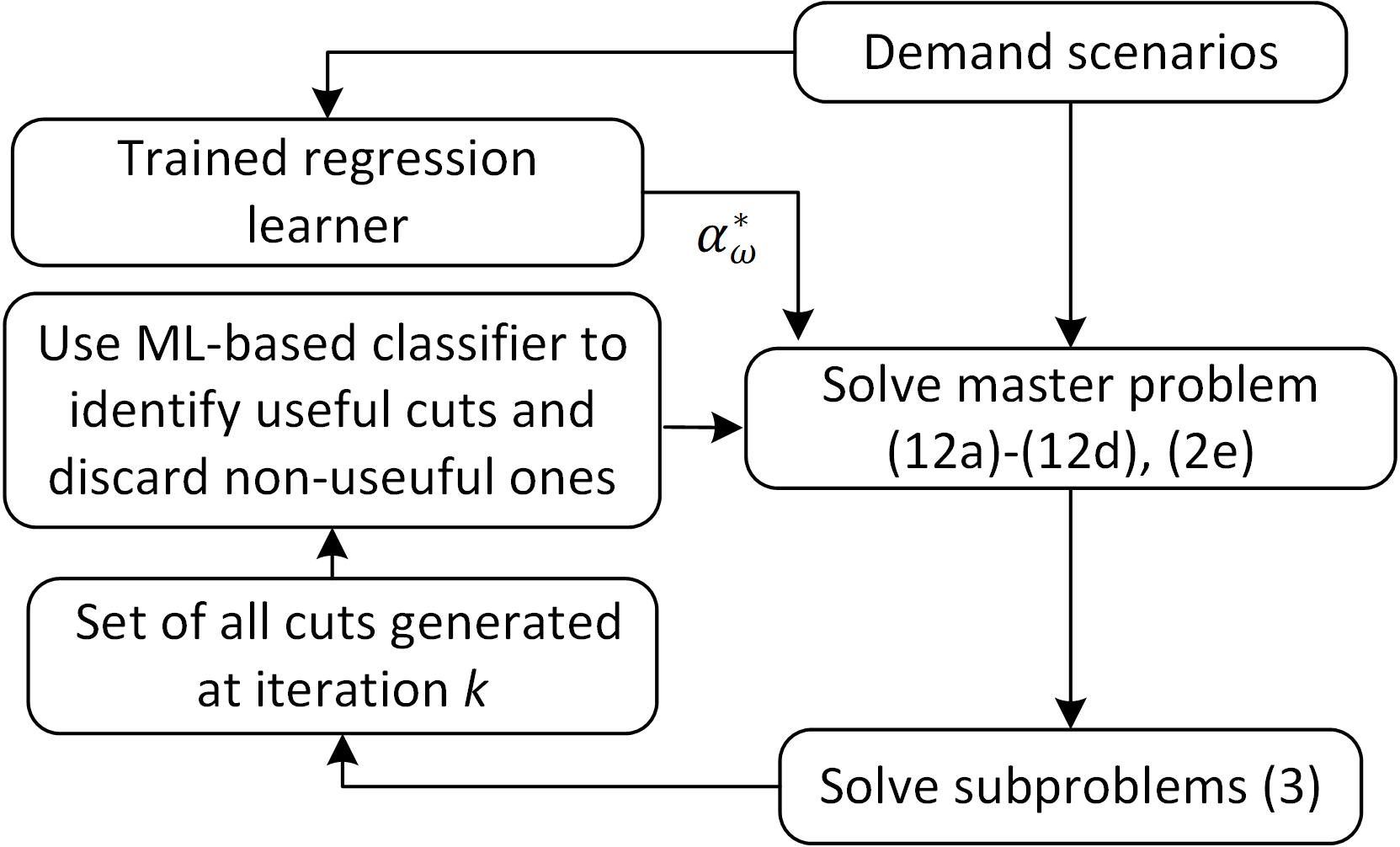}
\caption{CR-Benders flowchart.}
\label{fig8}
\end{figure}

\subsection{Subproblem Acceleration}
Two non-binding constraint removal strategies can be used to accelerate subproblem solutions, particularly for large SCUC problems. In the first learning-aided strategy, two classifiers can be trained to identify inactive (aka non-binding) network and generation ramp constraints at each Benders iteration and remove them from the model \cite{hasan2022topology}. Another more straightforward yet less efficient strategy is to solve a relaxed subproblem without network constraints. The relaxed subproblem solution and shift factors are used to check network constraints. The violated constraints are added to the relaxed subproblem, which is repeated until all network constraints are satisfied \cite{tejada2017security}. These two strategies reduce computation time and memory usage. We have used the second strategy.

\section{NUMERICAL RESULT}
The effectiveness of the proposed algorithms is evaluated in comparison to conventional multi-cut Benders decomposition. The YALMIP toolbox and CPLEX are used to implement Benders \cite{lofberg2004yalmip}, and PyTorch is used to construct neural network models. Simulations are conducted on a computer with an Intel(R) Xeon(R) 2.10 GHz CPU and 512 GB of RAM.

\subsection{Test Systems and Data Preparation}
Three test systems with various scheduling horizons are used, including the IEEE 24-bus, 118-bus, and 1354-bus systems \cite{babaeinejadsarookolaee2019power}. The 24-bus system includes ten generators and 34 transmission lines. The 118-bus system consists of 54 generators and 186 lines. The 1354-bus system has 270 generators and 1991 branches. Each case includes 40 equiprobable load profile scenarios. Table II shows the range of load perturbation with respect to base case demand values to incorporate daily and hourly load uncertainty. The convergence tolerance is set to 1\% for all cases.

\begin{table}
\setlength{\arrayrulewidth}{0.3mm}
\begin{center}
\caption{Parameter Range for Demand Dataset Generation}
\vspace{-5pt}
\label{table2}
\begin{tabular}{ m{2cm}  m{2cm}  m{2cm} } \hline
System & Load \\ \hline
 & [\(\eta_{s}^L\), \(\eta_{s}^U\)] & [\(\eta_{\omega}^L\), \(\eta_{\omega}^U\)] \\ \hline
24-bus & 70\%-130\% & 95\%-105\% \\ \hline
118-bus & 70\%-130\% & 95\%-105\% \\ \hline
1354-bus & 70\%-110\% & 95\%-105\% \\ \hline
\end{tabular}
\end{center}
\vspace{-15pt}
\end{table}

\subsection{Useful Cut Statistics}
Table III shows the average number of total cuts and useful cuts for different systems where $h$ stands for horizons. To determine the usefulness of a cut, we have calculated the contribution of each cut to improving the objective value of the master problem. A cut is useful if it results in a non-zero increase in the objective value. The percentage of the number of useful cuts to total cuts reduces as the size of the system and the number of scheduling horizons increase. For instance, for the IEEE 118-bus system with three time periods, only 13\% of cuts are labeled as useful, and the rest of 87\% are non-useful. Filtering non-useful cuts reduces the computational burden of the master problem significantly.

\begin{table}
\setlength{\arrayrulewidth}{0.3mm}
\begin{center}
\caption{Comparison of Total Cuts and Useful Cuts}
\vspace{-10pt}
\label{table3}
\begin{tabular}{ m{2cm}  m{2cm}  m{2cm} m{2cm} } \hline
Test case & Average number of total cuts & Average number of useful cuts & \% of useful cuts\\ \hline
24-bus,1h & 224 & 103 & 46\% \\ \hline
24-bus,12h & 1872 & 220 & 12\%  \\ \hline
118-bus,1h & 308 & 109 & 35\%  \\ \hline
118-bus,3h & 2280 & 288 & 13\%  \\ \hline
1354-bus,2h & 5520 & 795 & 15\%  \\ \hline
\end{tabular}
\end{center}
\vspace{-20pt}
\end{table}

\subsection{Master Problem Runtime and Solution Quality}
Table IV reports the average solver time of the master problem and the percentage of time improvement obtained by the three proposed approaches as compared to the conventional multi-cut Benders. During each iteration of the proposed approaches, we tackle a significantly smaller problem, resulting in considerable time savings. All three proposed approaches outperform conventional Benders. The time saving becomes more significant as the size of the optimization model increases. For instance, R-Benders reduces the master problem solution time from 236 seconds to 30 for the 118-bus system, an 87\% improvement. The time-saving of CR-Benders is more than R-Benders and C-Benders.

\begin{table}
\setlength{\arrayrulewidth}{0.3mm}
\begin{center}
\caption{Average Solver Time (and Improvement Percentage) Comparison in Seconds}
\vspace{-10pt}
\label{table4}
\begin{tabular}{ m{1.5cm}  m{1cm}  m{1.5cm} m{1.5cm} m{1.5cm} } \hline
System	& Benders	& R-Benders	& C-Benders & CR-Benders\\ \hline
24-bus,1h & 2.1	& 1.4 (33\%) & 0.73 (65\%) & 0.66 (69\%) \\ \hline
24-bus,12h & 37.7 &	32.9 (13\%)	& 9.7 (74\%) & 9.4 (75\%)  \\ \hline
118-bus,1h & 2.8 & 2.6 (7\%) & 2.5 (11\%) & 2.5 (11\%) \\ \hline
118-bus,3h & 236 & 30 (87\%) & 57 (76\%) & 24 (90\%) \\ \hline
1354-bus,2h & 110 & 74 (33\%) & 65 (41\%) & 62 (44\%)  \\ \hline
\end{tabular}
\end{center}
\end{table}

The average number of Benders iterations is shown in Fig. 9. The three proposed approaches, particularly C-Benders and CR-Benders, take roughly the same number of iterations as conventional Benders decomposition. A comparison of the solution time and the number of iterations, as well as the number of useful cuts from Table III, show that R-Benders, C-Benders, and CR-Benders can capture the necessary information to build a computationally less expensive master problem.

\begin{figure}
\centering
    \includegraphics[width=.3\textwidth]{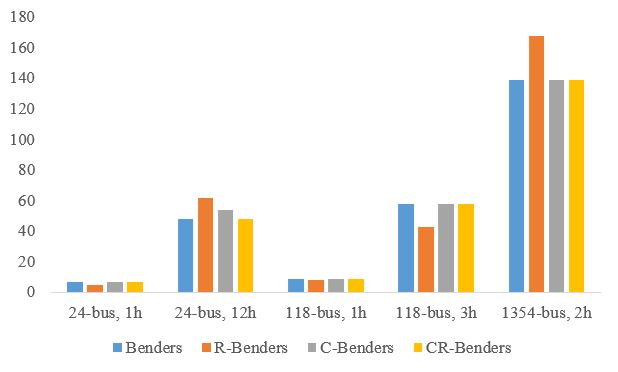}
\caption{Average number of iterations.}
\label{fig9}
\vspace{-15pt}
\end{figure}

We use a CostGap index to assess the solution quality. The CostGap index measures the difference between the objective values generated by the proposed approaches $f^p$ and the conventional Benders $f^BD$.

\begin{equation} \tag{13}
    \label{13}
    CostGap\% = \frac{|f^P-f^{BD}|}{f_{BD}} \times  100
\end{equation}

Table V presents the average CostGap index for all cases. The average cost gap for all test cases is negligible. For instance, the gap is less than 0.02\% for the 118-bus system with a 3-hour horizon. We have observed several cases for which the proposed approaches obtain even a better solution than the conventional Benders. This is due to a warm start with stronger cuts in the first iteration.

\begin{table}
\setlength{\arrayrulewidth}{0.3mm}
\begin{center}
\caption{Costgap Index (\%)}
\vspace{-10pt}
\label{table5}
\begin{tabular}{ m{1.5cm}  m{1.3cm}  m{1.3cm}  m{1.4cm} } \hline
System	& R-Benders	& C-Benders	& CR-Benders\\ \hline
24-bus,1h & 0 & 0.13	& 0.07 \\ \hline
24-bus,12h & 0.05 & 0.8 &	0.03 \\ \hline
118-bus,1h & 0.6 & 0.01	& 0.01\\ \hline
118-bus,3h & 0 & 0.02 & 0.02 \\ \hline
1354-bus,2h & 0.04 & 0.05 & 0 \\ \hline
\end{tabular}
\end{center}
\vspace{-20pt}
\end{table}

\subsection{Memory Usage }

\begin{figure}
\centering
    \includegraphics[width=.3\textwidth]{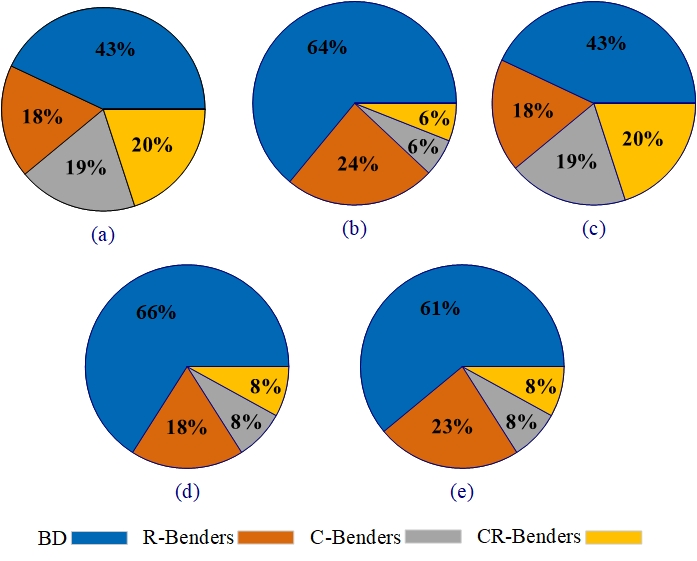}
\caption{Memory usage comparison. A) 24-bus (1h), b) 24-bus (12h), c) 118-bus (1h), d) 118-bus (3h), and e) 1354-bus (2h).}
\label{fig10}
\end{figure}

The memory usage of the proposed approach is compared with the conventional Benders. The memory occupied to form the constraint set from Benders cuts is illustrated in Fig. 10. The memory requirement reduces significantly by filtering out non-useful cuts from the model. For example, Fig. 10.d shows that for the 118-bus system with 3h horizon, the memory usage of CR-Benders is 87.5\% less than that of the conventional Benders decomposition. Fig. 11 demonstrates the number of cuts added to the master problem per iteration. The slope of the increment for the conventional Benders is steeper than that of RC-Benders.

\begin{figure}
\centering
    \includegraphics[width=.3\textwidth]{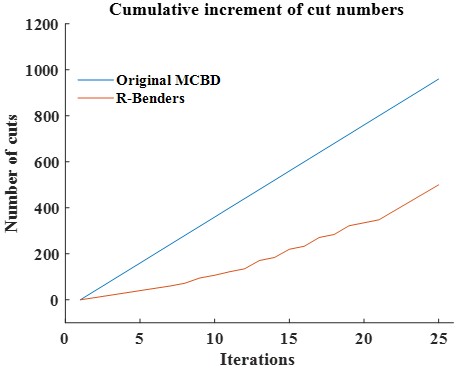}
\caption{Cumulative increment of the number of cuts for IEEE 118-bus 3h case.}
\label{fig11}
\vspace{-15pt}
\end{figure}
\vspace{-9pt}
\subsection{Non-convergent case analysis}
In specific situations, the conventional Benders decomposition (BD) technique struggles to achieve convergence within the allocated timeframe, even in the presence of an optimal solution. Through our experiments, we have demonstrated that the utilization of a predicted subproblem proxy can effectively address these scenarios and enhance the duality gap. We focused on two particular scenarios involving a 118-bus system and a 3-hour horizon, where the original MCBD approach failed to converge within 400 iterations. By employing the R-benders method, we observed an average convergence within 59 iterations for these same scenarios. It is evident that the incorporation of the proxy value yields superior performance.
\vspace{-9pt}
\section{CONCLUSION}
The proposed algorithms aim to reduce the computational costs of two-stage SCUC using machine learning-based warm-start and cut filtering. The proposed algorithms provide guaranteed feasible and near-optimal solutions. Conventional Benders decomposition suffers from slow convergence and might not reach an optimal point within a specified time. Not all cuts generated through Benders iterations contain useful information to form the feasible region of the master problem. Removing non-useful cuts and forming stronger cuts enhance Benders decomposition performance. We trained a regressor to predict subproblem objective proxy variables for the master problem. The predicted values form the first set of cuts in the first Benders iteration. After solving the first iteration and onward, cuts are filtered by comparing the numerical value of cuts and the numerical value of proxy variables obtained from that particular iteration. A classifier is trained to identify useful cuts automatically after each Benders iteration. The proposed approaches are applicable to many MIP problems.
Simulation results on various test systems show CR-Benders outperform R-Benders and C-Benders in terms of solution time and memory saving. On average, CR-Benders leads to 58\% time saving and 77\% memory saving, which is higher than the other two approaches. While C-Benders and CR-Benders work better than R-Benders, they need a computationally expensive dataset generation for classification learning training. This procedure is, however, offline and needs to be implemented once.

\bibliographystyle{ieeetr}
\bibliography{BendersDecomposition.bib}

\end{document}